\documentclass[12pt]{article}
\usepackage{amssymb}
\usepackage{amsmath}
\usepackage{cases}
\usepackage{version}
\usepackage{graphicx}

\newtheorem{lem}{Lemma}[section]
\newtheorem{Prop}{Proposition}[section]
\newtheorem{Theor}{Theorem}[section]

\newtheorem{rmrk}{Remark}[section]

\newcommand{\cqfd}{\hfill $\square$}
\newcommand{\R}{\mathbb R}

\def \E{\mathbb E}

\title{About the Stein equation for the generalized inverse Gaussian and Kummer distributions}
\author{Essomanda Konzou, \ Angelo Efo\'evi Koudou   \\
 {\it \small Institut Elie Cartan de Lorraine, UMR CNRS 7502,  Universit\'e de Lorraine,  Universit\'e de Lom\'e }} 
\date{July 31, 2018}
\begin{document}
\maketitle

%
%
\begin{abstract} We propose a Stein characterization of the Kummer distribution on $(0,\infty)$. This result follows from our observation  that the density of the Kummer distribution satisfies a certain differential equation, leading to a solution of the related Stein equation. A bound is derived  for the solution, under a condition on the parameters. The derivation of this bound is carried out  using the same framework  as in Gaunt 2017 [A Stein characterisation of the generalized hyperbolic distribution. {\it ESAIM: Probability and Statistics}, {\bf 21}, 303--316]  in the case of the generalized inverse Gaussian distribution, which we revisit by correcting a minor error in the latter paper. 
 \end{abstract}
%
%

\section{Introduction}

For $a >0$, $b\in \R$, $c >0$, the Kummer distribution with parameters $a, b, c$ has density 
$$ k_{a, b, c}(x)=
\frac{ 1 }{\Gamma (a)\psi (a, a-b+1;c)}  x^{a-1}   (1+x)^{-a-b} e^{- c x}, \ (x>0) $$
where $\psi$ is the confluent hypergeometric function of the second kind. 

The generalized inverse Gaussian (hereafter $GIG$) distribution with parameters $p\in\mathbb{R}$, $a > 0$, $b > 0$ has density
 $$
g_{p,a, b}(x)= \dfrac{\left( a/b\right)^{p/2}}{2K_{p}(\sqrt{ab})}x^{p-1}e^{-\frac{1}{2}\left( ax+b/x\right)}, \quad x>0,
$$
where $K_p$ is the modified Bessel function of the third kind. 

For details on GIG and Kummer distributions see \cite{Koudou caraterisation, kv, Weso} and references therein, where one can see for instance that these distributions are involved in some characterization problems related to the so-called Matsumoto-Yor property. 

In this paper, these two distributions are considered  in the context of Stein's method.
This method introduced in \cite{Stein} is a technique used to bound the error in the approximation of the distribution of a random variable of interest by another probability (for instance the normal) distribution. For an overview of Stein's method see \cite{Chen, Ross}.  The first steps of this method consist in finding an operator called Stein operator characterizing the targeted distribution, then solving  the corresponding so-called {\it Stein equation}. 
 
 One finds in  \cite{Stein} a seminal  instance of the method, where Stein showed that a random variable $X$  has a standard normal distribution  if and only if for all real-valued absolutely continuous function $f$ such that $\E \left| f^\prime(Z)\right| <\infty$ for $Z\sim N(0,1)$,
 $$\E\left[ f^\prime(X)-Xf(X)\right] =0.$$
 The corresponding {\it  Stein equation}  is 
 $$ f^\prime(x)-xf(x)=h(x)-\E h(Z)$$
 where $h$ is a bounded function and $Z$ a random variable  following the standard normal distribution.
The operator $f\mapsto T_f$  defined by $(T_f)(x)=f^\prime(x)-xf(x)$ is the corresponding {\it Stein operator}.

 If a function $f_h$ is a solution of the previous equation, then for any random variable $U$ we have 
$$ |\E (f_h^\prime(U)-Uf_h(U))|=|\E (h(U))-\E (h(Z))|.$$ Thus, in order to bound $|\E (h(U))-\E (h(Z))|$ given $h$, its enough to find a solution $f_h$ of the Stein equation and to bound the left-hand side of the previous equation. 
The problem of solving the Stein equation for other distributions than the standard normal distribution and  bounding the solution and its derivatives has been widely studied in the literature (see \cite{Goldstein} among many others).

 The aim of this paper is to solve the Stein equation and derive a bound of the solution for the Kummer distribution (which is new) and for the generalized inverse Gaussian distribution (which has been done in \cite{Gaunt}, but there was a little mistake in the bound of the solution).

The idea of this paper emerged  by reading the remarkable work by  \cite{Gaunt} about a Stein characterization of the generalized hyperbolic distribution of which the generalized inverse Gaussian distribution (GIG)  is a limiting case. Among many other results,  \cite{Gaunt}  solved the GIG Stein equation and bounded the solution by using a general result obtained in   \cite{Schoutens 2001} when the targeted distribution has a density $g$ satisfying 
\begin{equation}
\label{defstau}
\left( s (x)g(x)\right) ^\prime=\tau (x)g(x)
\end{equation}
for some polynomial functions $s$ and $\tau$.  Also a bound was obtained for the solution under the condition that the function $\tau$ be a decreasing linear function. 
But since this linearity condition does not hold in the GIG case, the bound given by \cite{Gaunt} has to be slightly corrected. This is done in Theorem \ref{solutionGIG} after recalling the general framework of Schoutens \cite{Schoutens 2001} and adapting it to the cases where $\tau$  is decreasing but not necessarily linear. Indeed, we realized that the procedure adopted in \cite{Schoutens 2001} still works, via a slight change, even if $\tau$ is not linear.  

Observing that the Kummer density also satisfies  \eqref{defstau}, we can use the same methodology (Theorem \ref{theoKummer})  for this distribution.
 We have to put the restrictions $p\leq -1$ for the GIG density and  $1-b-c\leq 0$ for the Kummer density in order for the corresponding function $\tau$ to be decreasing on $(0, \infty)$.

In Section 1 we recall the general framework established by \cite{Schoutens 2001} for densities $g$ satisfying  \eqref{defstau} without the assumption of linearity of $\tau$. We retrieve the Stein operator given in  \cite{Schoutens 2001} by using the density approach  initiated in \cite{Steinetal} and further developed in \cite{LeySwan}.

In Section 2 we show the application of this method to the GIG distribution as mentioned in \cite{Gaunt} by giving the right bound for the solution of the Stein equation.
Section 3 is devoted to the Stein characterization and the Stein equation related to the Kummer distribution.

\section{Stein characterization in the Schoutens framework}
\setcounter{equation}{0}
Theorem 1 in \cite{Schoutens 2001} addressed the problem of establishing a Stein characterization for probability distributions with density $g$ satisfying
\eqref{defstau}
for some polynomial functions $s$ et $\tau$, and proved that a Stein operator in this case is  $f\mapsto sf^\prime+ \tau f$. We realized (see the following theorem) that the same Stein operator can be arrived at by using the density approach designed in \cite{Steinetal} and \cite{LeySwan}. The support of the density may be any interval, but here we take this support to be $(0,\infty)$ in the purpose of the application to the  GIG and  Kummer distributions.

\begin{Theor} \label{steingen}
Consider a density $g$ on $(0,\infty)$ such that \eqref{defstau} holds for some polynomial functions $s$ and $\tau$. 
Then a positive random variable $X$  has density $g$ if and only if for any differentiable  function $f$ such that 
$\lim\limits_{\substack{x\to 0}} s(x)g(x)f(x)=\lim\limits_{\substack{x\to \infty}} s(x)g(x)f(x)=0$,
$$
\E \left[s(X)f^\prime (X) + \tau (X) f(X)  \right]=0.
$$
\end{Theor}

{\bf Proof:}
We use Corollary 2.1 of \cite{LeySwan}. According to this corollary, a Stein operator  related to the density function $g$ is 
$$T _g f (x)= \frac{1}{g(x)} (fg)^\prime (x.) $$
Applyng this operator to $sf$, we have 
\begin{eqnarray}
  T _g (sf) (x)& = &\frac{1}{g(x)} (sfg)^\prime (x) \nonumber \\
& = & \frac{1}{g(x)}(f^\prime (x)s(x)g(x) +f(x)  (sg)^\prime (x)) \nonumber
\end{eqnarray}
which, by \eqref{defstau}, reads
\begin{eqnarray}
  T_g (sf) (x)& = & \frac{1}{g(x)}(f^\prime (x)s(x)g(x) +f(x)  \tau(x) g(x)) \nonumber \\
& = &  f^\prime (x)s(x) +f(x)  \tau(x). \nonumber \end{eqnarray}

\cqfd

Theorem \ref{steingen} shows that  the Stein equation related to any density $g$ satisfying  \eqref{defstau} enjoys the tractable form 
\begin{equation}
\label{stein}
s(x)f^\prime(x)+\tau(x)f(x)=h(x)-\E h(W)
\end{equation}
where $W$ is random variable with density $g$.
Schoutens \cite{Schoutens 2001} found a solution to the Stein equation (\ref{stein})  and established  a bound for the solution, under the condition that the function $\tau$ be a decreasing linear function (which is the case for the so-called  Pearson and Ord classes of distributions considered in \cite{Schoutens 2001}). 

The following result comes from Proposition 1 in \cite{Schoutens 2001}.
We  again take the support of the density function to be $(0,\infty)$.

\begin{Prop}
		\label{caracterisation generale prem}	
Consider  a density function $g>0$ on $(0,\infty)$ satisfying Equation (\ref{defstau}), for some  polynomial functions  $s$ and $\tau$.
Then a solution of the Stein equation \eqref{stein}   is
	\begin{equation}
	\label{solution gen}
	\begin{split}
	f_h(x)= & \dfrac{1}{ s(x)g(x)}\displaystyle\int_{0}^{x}g(t)\left[ h(t)-\E h(W)\right] dt\\
	= & \dfrac{-1}{ s(x)g(x)}\displaystyle\int_{x}^{+\infty}g(t)\left[ h(t)-\E h(W)\right] dt.\\
	\end{split}
	\end{equation}
\end{Prop}

\begin{rmrk}
The proof  of this proposition is given in \cite{Schoutens 2001} just by calculating the derivative of the function $f_h$ defined  by \eqref{solution gen} and checking  directly that $f_h$ satisfies \eqref{stein}. Our following proposition complements this result.
\end{rmrk}

\begin{Prop}
		\label{caracterisation generale}	
Under the notation and assumptions of Proposition \ref{caracterisation generale prem},
\begin{itemize}
\item
The  solutions of  the Stein equation \eqref{stein}  are of the form
	\begin{equation}
	\label{solution gen bis}
	\begin{split}
	f_h(x)=  & \dfrac{1}{ s(x)g(x)}\displaystyle\int_{0}^{x}g(t)\left[ h(t)-\E h(W)\right] dt   +\dfrac{C}{ s(x)g(x)}\\
	= & \dfrac{-1}{ s(x)g(x)}\displaystyle\int_{x}^{+\infty}g(t)\left[ h(t)-\E h(W)\right] dt +\dfrac{C}{ s(x)g(x)}\\
	\end{split}
	\end{equation}
where $C$ is constant.
\item Suppose $\lim\limits_{\substack{x\to 0}} s(x)g(x)=0$. For the solution to be bounded, it is necessary that $C=0$ in \eqref{solution gen bis}.
\end{itemize}
\end{Prop}

{\bf Proof:}


Multiplying both sides of \eqref{stein} by $g(x)$ we have  

$$s(x)g(x)f^\prime(x)+\tau(x)g(x)f(x)=g(x) (h(x)-\E h(W))$$
which, by \eqref{defstau}, can be written
$$s(x)g(x)f^\prime(x)+(sg)^\prime(x)f(x)=g(x) (h(x)-\E h(W)),$$
i.e.
$$(sgf)^\prime(x)=g(x) (h(x)-\E h(W)).$$
As a consequence, there exists a constant $C$ such that
\begin{equation}
	\label{const}
s(x)g(x)f(x)= \displaystyle\int_{0}^{x}g(t)\left[ h(t)-\E h(W)\right] dt   + C \end{equation}
which implies \eqref{solution gen bis}.

Suppose $f$ is bounded. Since $\lim\limits_{\substack{x\to 0}} s(x)g(x)=0$, letting $x$ tend to $0$ in \eqref{const} yields $C=0$.


The second expression for $f_h$ follows from the fact that, since $W$ has density $g$,
$$
 \displaystyle\int_{0}^{+\infty}g(t)\left[ h(t)-\E h(W)\right] dt =0.$$ 
\cqfd

The following proposition proves that the solution given by  \eqref{solution gen}  is bounded  indeed if $h$ is bounded, and thus is the unique bounded solution to the Stein equation associated to the density $g$.   A  bound is provided. 

\begin{Prop} 
\label{borne}	
Consider  a density function $g>0$ on $(0,\infty)$ satisfying Equation (\ref{defstau}), where $s$ and $\tau$ are  polynomial functions such that  $s>0$ on  $(0,\infty)$ and $\tau$ is decreasing and has a unique zero $\alpha$  on $(0,\infty)$. 
Assume that $\lim\limits_{\substack{x\to 0}} s(x)g(x)=\lim\limits_{\substack{x\to \infty}} s(x)g(x)=0$.
	If $h$ is a bounded continuous function, then
	\begin{equation}
	\label{bornefh}
	\left| \left| f_h\right| \right| \leq M\left| \left| h(.)-\E h(Z)\right| \right| 
	\end{equation}
	where 
	$$ M=\max \left( \dfrac{1}{s(\alpha)g(\alpha)} \displaystyle\int_{0}^{\alpha}g(t)dt ;
	\dfrac{1}{s(\alpha)g(\alpha)} \displaystyle\int_{\alpha}^{+\infty}g(t)dt\right)  $$
	and 
	$ \left| \left| f_h\right| \right|=\displaystyle\sup\limits_{x>0}\left| f_h(x)\right|.$
\end{Prop}

\begin{rmrk}
This result is a reformulation of  Lemma 1 in  \cite{Schoutens 2001} without  the assumption that $\tau$ is linear. With this assumption,\cite{Schoutens 2001} established the same bound with $\alpha = \E (X)$ (for a random variable $X$ with density $g$),   which is not true if $\tau$ is not linear.
The proof given below follows the lines of that of \cite{Schoutens 2001} where we observed that the assumption of linearity of $\tau$ was used nowhere except to state that its only zero is $\alpha = \E(X)$.
\end{rmrk}

The proof of Proposition \ref{borne} uses the following lemma :

\begin{lem}
	\label{lemme}
	Under the assumptions of Proposition \ref{borne},
	$$ \displaystyle\int_{0}^{x}g(t)dt\leq \dfrac{s(x) g(x)}{\tau(x)}\quad \text{ for } \quad x<\alpha$$
	and
$$ \displaystyle\int_{x}^{+\infty}g(t)dt\leq \dfrac{-s(x) g(x)}{\tau(x)} \quad \text{ for }\quad  x>\alpha.$$\\	
\end{lem}

{\bf Proof:}
Suppose $x<\alpha$. Since  $\tau$ is positive and decreasing on $(0,\alpha )$, we have
$ \dfrac{\tau (t)}{\tau (x)}\geq 1$ for all $ t\leq x $.
Therefore 
\[
\begin{split}
\displaystyle\int_{0}^{x}g(t)dt \leq & \displaystyle\int_{0}^{x}g(t)\dfrac{\tau (t)}{\tau (x)}dt\\
 = & \dfrac{1}{\tau (x)}\displaystyle\int_{0}^{x}\tau(t)g(t)dt\\
 = & \dfrac{s(x) g(x)}{\tau(x)}
\end{split}\]
because of \eqref{defstau} and as $\lim\limits_{\substack{t\to 0}} s(t)g(t)=0$.

For $x>\alpha$,  since $\tau$ is negative and decreasing on $(\alpha, \infty )$, we have $\dfrac{\tau (t)}{\tau (x)}\geq 1$  for all
$ t\geq x$. As a consequence, 
\[
\begin{split}
\displaystyle\int_{x}^{+\infty}g(t)dt \leq & \displaystyle\int_{x}^{+\infty}g(t)\dfrac{\tau (t)}{\tau (x)}dt\\
= & \dfrac{-s(x) g(x)}{\tau(x)}
\end{split}\] since $\lim\limits_{\substack{t\to \infty}} s(t)g(t)=0$. \cqfd

Now, let us prove \eqref{bornefh}.

{\bf Proof:}

For $x<\alpha$,

\[
\begin{split}
\left| f_h(x)\right|= & \left| \dfrac{1}{ s(x)g(x)}\displaystyle\int_{0}^{x}g(t)\left[ h(t)-\E h(W)\right] dt\right| \\
                 \leq & \dfrac{1}{ s(x)g(x)}\displaystyle\int_{0}^{x}g(t)\left| h(t)-\E h(W)\right| dt\\
                 \leq & \left|\left|  h(.)-\E h(W)\right| \right|\dfrac{1}{ s(x)g(x)}\displaystyle\int_{0}^{x}g(t) dt.
\end{split}
\]
Let $l(x)=\dfrac{1}{s(x)g(x)}\displaystyle\int_{0}^{x}g(t) dt$. Then $l$ is differentiable on $(0,\infty)$ and
\[
\begin{split}
l^\prime(x)= & \dfrac{-\left( s(x)g(x)\right)^\prime}{\left(  s(x)g(x)\right)^2}\displaystyle\int_{0}^{x}g(t) dt+
\dfrac{1}{ s(x)}\\
  = & \dfrac{- \tau(x)}{  s^2(x)g(x)}\displaystyle\int_{0}^{x}g(t) dt+ \dfrac{1}{ s(x)}
\end{split}
\]
Using Lemma \ref{lemme}, we conclude that $l^\prime(x)\geq 0$. Then $l(x)\leq l(\alpha).$

For $x>\alpha$,

\[
\left| f_h(x)\right|=  \left|\left|  h(x)-\E h(W)\right| \right|\dfrac{1}{ s(x)g(x)}\displaystyle\int_{x}^{+\infty}g(t) dt.
\]
Let $u(x)=\dfrac{1}{ s(x)g(x)}\displaystyle\int_{x}^{+\infty}g(t) dt$. The function $u$ is differentiable on  $(0,\infty)$ and
\[
u^\prime(x)= \dfrac{- \tau(x)}{  s^2(x)g(x)}\displaystyle\int_{x}^{+\infty}k(t) dt-\dfrac{1}{ s (x)}
\]
By Lemma \ref{lemme}, we conclude that $u^\prime(x)\leq 0$. Then $u(x)\leq u(\alpha).$ \cqfd

In the two next sections we apply the previous results to the GIG and Kummer distributions.

\section{About the Stein equation of the generalized inverse Gaussian distribution}
\setcounter{equation}{0}

Recall that the density of the GIG distribution with parameters $p\in\mathbb{R}$, $a > 0$, $b > 0$  is  
	\begin{equation}
\label{dgig2}
g_{p,a, b}(x)= \dfrac{( a/b)^{p/2}}{2K_{p}(\sqrt{ab})}x^{p-1}e^{-\frac{1}{2}\left( ax+b/x\right)}, \quad x>0,
\end{equation}
where $K_p$ is the modified Bessel function of the third kind. 

Let 
\begin{equation}
\label{gig tau}
 s(x)=x^2 \quad \text{ and }\quad \tau (x)=\frac{b}{2}+(p+1)x-\frac{a}{2}x^2.
\end{equation}
Then, as observed by \cite{Gaunt},  the GIG density $g_{p,a, b}$  satisfies $$\left( s (x)g_{p,a, b}(x)\right) ^\prime=\tau (x)g_{p,a, b}(x).$$
This enables us to apply  Theorem \ref{steingen} to retrieve the following  Stein characterization of the GIG distribution given in \cite{Koudou caraterisation} and \cite{Gaunt}:  


\begin{Prop}
		\label{caracterisation gig}	
	A random variable $X$ follows the $GIG$ distribution with density $g_{p,a, b}$  if and only if, for all real-valued and differentiable function $f$ such that\\
		$\lim\limits_{\substack{x\to\infty}}g_{p,a, b}(x)f(x)=\lim\limits_{\substack{x\to 0}}g_{p,a, b}(x)f(x)=0,$ we have:
$$		\mathbb{E}\left[ X^2f^\prime (X)+\left( \frac{b}{2}+(p+1)X-\frac{a}{2}X^2\right) f(X)\right] =0.$$	
\end{Prop}

The corresponding  Stein equation is 
\begin{equation}
\label{eq diff gig}
x^2f^\prime (x)+\left( \frac{b}{2}+(p+1)x-\frac{a}{2}x^2\right) f(x)= h(x)-\E h(W)
\end{equation}
 where $h$ is  a bounded function and $W$ a random variable  following the $GIG$ distribution with parameters $p,a, b$.\\

We apply Proposition \ref{caracterisation generale} and Proposition \ref{borne} to solve Equation \eqref{eq diff gig} and bound the solution. Let us check that the assumptions of these propositions are true in the GIG case.

Firstly, we note  that, by \eqref{dgig2},  $$s(x) g_{p,a, b}(x)= \dfrac{\left( a/b\right)^{p/2}}{2K_{p}(\sqrt{ab})}x^{p+1}e^{-\frac{1}{2}\left( ax+b/x\right)}, \quad x>0,$$
which shows that $\lim\limits_{\substack{x\to\infty}}s(x)g_{p,a, b}(x)=\lim\limits_{\substack{x\to 0}}s(x)g_{p,a, b}(x)=0$.

Secondly, observe  that if $p\leq -1$, then the function $\tau$ defined by \eqref{gig tau} is decreasing on $(0,\infty)$ and that its only zero on $(0,\infty)$  is $\alpha=\dfrac{p+1+\sqrt{(p+1)^2+ab}}{a}.$ 

Thus, by using Proposition \ref{caracterisation generale} and Proposition \ref{borne}, we obtain the following theorem.

\begin{Theor} \label{solutionGIG}
	The GIG Stein equation \eqref{eq diff gig}  has solution 
	\begin{equation} \label{fh}
	\begin{split}
	f_h(x)= & \dfrac{1}{ s(x)g_{p,a, b}(x)}\displaystyle\int_{0}^{x} g_{p,a, b}(t)\left[ h(t)-\E h(W)\right] dt\\
	= & \dfrac{-1}{  s(x)g_{p,a, b}(x)}\displaystyle\int_{x}^{+\infty}g_{p,a, b}(t)\left[ h(t)-\E h(W)\right] dt\\
	\end{split}
	\end{equation} where $W$ follows the $GIG$ distribution with parameters $p,a, b$.

	If $h$ is a bounded continuous function and $p\leq -1$, then the function defined by \eqref{fh} is the unique bounded solution of  \eqref{eq diff gig} and
$$	\left| \left| f_h\right| \right| \leq M\left| \left| h(.)-\E h(W)\right| \right| $$
	where 
$$\alpha=\dfrac{p+1+\sqrt{(p+1)^2+ab}}{a},$$
	$$ M=\max \left( \dfrac{1}{s(\alpha)g_{p,a, b}(\alpha)} \displaystyle\int_{0}^{\alpha}g_{p,a, b}(t)dt ;
	\dfrac{1}{s(\alpha)g_{p,a, b}(\alpha)} \displaystyle\int_{\alpha}^{+\infty}g_{p,a, b}(t)dt\right).  $$
\end{Theor}

\begin{rmrk}
This result was claimed by Gaunt (see \cite{Gaunt}) with $\alpha=\E (X)$ by applying Proposition $1$ of   \cite{Schoutens 2001}. The only slight mistake  is that $\tau$ is not a  polynomial function of degree one as in \cite{Schoutens 2001}.
\end{rmrk}

\section{About the Stein equation related to the Kummer distribution}
\setcounter{equation}{0}
Recall that for $a >0$, $b\in \R$, $c >0$, the Kummer distribution $K(a, b, c)$ has density $$ k_{ a, b, c}(x)=
\frac{ 1 }{\Gamma (a)\psi (a, a-b+1;c)}  x^{a-1}   (1+x)^{-a-b} e^{- c x}, \ (x>0)$$
where $\psi$ is the confluent hypergeometric function of second kind. 
Let 
\begin{equation}
\label{lu}
s(x)=x(1+x) \quad \text{ and }\quad \tau (x)= (1-b)x-cx(1+x)+a.
\end{equation}
 Then
$$\left( s (x)k_{ a, b, c}(x)\right) ^\prime=\tau (x)k_{ a, b, c}(x.)$$
Then we can use  Theorem \ref{steingen} to obtain the following Stein characterization of the Kummer distribution:
\begin{Theor}
\label{caracterisation 1}	
A random variable $X$ follows the $K(p,a,b)$ distribution if and only if, for all differentiable function $f$,
$$\E \left[X(X+1)(f^\prime (X) + \left( (1-b)X-cX(1+X)+a   \right) f(X)  \right]=0.$$
\end{Theor}

The  corresponding Stein equation  is 
\begin{equation}
\label{eqdiff 2}
x(x+1)f^\prime (x) + \left[  (1-b)x-cx(1+x)+a  \right]  f(x)=h(x)-\E h(W)
\end{equation} 
where $W$ has density $k_{ a, b, c}$.

We have $$s(x)k_{ a, b, c}(x)=\frac{ 1 }{\Gamma (a)\psi (a, a-b+1;c)}  x^{a}   (1+x)^{1-a-b} e^{- c x}, \  \ x>0$$
and we see that $\lim\limits_{\substack{x\to\infty}}s(x)k_{ a, b, c}(x)=\lim\limits_{\substack{x\to 0}}s(x)k_{ a, b, c}(x)=0$.

Note also that if $1-b-c\leq 0$, then the function $\tau$ defined by \eqref{lu} is decreasing on $(0,\infty)$ and its only zero on this interval  is $\alpha=\dfrac{1-b-c+\sqrt{(1-b-c)^2+4ac}}{2c}.$

Then we use again  Proposition \ref{caracterisation generale} and Proposition \ref{borne} to obtain the following result:

\begin{Theor} \label{theoKummer}
	The Kummer Stein equation \eqref{eqdiff 2}  has solution 
	\begin{equation} \label{fhk}
	\begin{split}
	f_h(x)= & \dfrac{1}{ s(x)k_{ a, b, c}(x)}\displaystyle\int_{0}^{x}k_{ a, b, c}(t)\left[ h(t)-\E h(W)\right] dt\\
	= & \dfrac{-1}{ s(x)k_{ a, b, c}(x)}\displaystyle\int_{x}^{+\infty}k_{ a, b, c}(t)\left[ h(t)-\E h(W)\right] dt\\
	\end{split}
	\end{equation} where $W\sim K(a,b,c)$.

	If $h$ is a bounded continuous function and $1-b-c\leq 0$, then $f_h$ defined by \eqref{fhk} is the unique bounded solution of \eqref{eqdiff 2} and
	$$\left| \left| f_h\right| \right| \leq M\left| \left| h(.)-\E h(W)\right| \right| $$
	where 
$$\alpha=\dfrac{1-b-c+\sqrt{(1-b-c)^2+4ac}}{2c},$$
	$$ M=\max \left( \dfrac{1}{s(\alpha)k_{ a, b, c}(\alpha)} \displaystyle\int_{0}^{\alpha}k_{ a, b, c}(t)dt ;
	\dfrac{1}{s(\alpha)k_{ a, b, c}(\alpha)} \displaystyle\int_{\alpha}^{+\infty}k_{ a, b, c}(t)dt\right).  $$
\end{Theor}

 \begin{rmrk} 
These results could be used in future work to provide rates of convergence in limit problems related to the GIG and Kummer distributions.
\end{rmrk}


\begin{thebibliography}{}

\end{thebibliography}


\begin{thebibliography}{99}
	

\bibitem{Chen}
 L. H. Y. Chen, L. Goldstein and  Q.-M. Shao.  {\it Normal approximation by Stein's method}. Probability and its Applications (New York).  Springer, Heidelberg, 2011.


 \bibitem{Gaunt}
 R. E. Gaunt. A Stein characterisation of the generalized hyperbolic distribution. {\it ESAIM: Probability and Statistics} {\bf 21} (2017) 303--316.


\bibitem{Goldstein} 
L. Goldstein and G. Reinert. Stein's method for the Beta distribution and the P\`{o}lya-Eggenberger urn. \emph{ Adv. Appl. Probab.} {\bf 50} (2013) 1187--1205.


\bibitem{Koudou caraterisation} 
A. E. Koudou and C. Ley.  Characterizations of GIG laws: a survey complemented with two new results \emph{Proba. Surv.} \textbf{11} (2014) 161--176.

\bibitem{kv}
A. E. Koudou  and P. Vallois.  Independence properties of the Matsumoto-Yor type. {\it  Bernoulli} {\bf 18}  (2012) 119--136.


\bibitem{LeySwan}
C. Ley  and Y. Swan.  Stein's density approach and information inequalities. \emph{Electron. Comm. Probab.} \textbf{18(7)} (2013) 1--14.


\bibitem{Ross}
N. Ross.  Fundamentals of Stein's method. \emph{Probab. Surv.} \textbf{8} (2011) 210--293.

\bibitem{Schoutens 2001} 
W. Schoutens.  Orthogonal Polynomials in Steins Method. \emph{J. Math. Anal. Appl.} \textbf{253} (2001) 515--531.


\bibitem{Stein}
 C.  Stein.  A bound for the error in the normal approximation to the distribution of a sum of dependent random variables. In \emph{Proceedings of the Sixth Berkeley Symposium on Mathematical Statistics and Probability} (Vol 2, pp. 583--602)  (1972). Berkeley: University of California Press.

\bibitem{Steinetal}
C. Stein, P. Diaconis, S. Holmes and  G. Reinert.  Use of exchangeable pairs in the analysis of simulations. In \emph{Persi Diaconis and Susan Holmes, editors, Stein's method: expository lectures and applications}, vol. 46 of \emph{IMS Lecture Notes Monogr. Ser},  1--26  (2004). Beachwood,
Ohio, USA: Institute of Mathematical Statistics.

\bibitem{Weso}
A. Piliszek and J. Weso\l owski. Change of measure technique in characterizations of the gamma and Kummer distributions. 
\emph{J. Math. Anal. Appl.} \textbf{458} (2018)  967--979.

\end{thebibliography}
\end{document}